\documentclass{article}
\usepackage{amsfonts}
\usepackage[latin5]{inputenc}
\usepackage{amssymb}
\usepackage{amsmath}

\numberwithin{equation}{section}
\input{tcilatex}
\begin{document}

\title{Special Smarandache Curves According To Darboux Frame In $E^{3}$}
\author{ÖZCAN BEKTAŞ $^{\ast }$ \and SALİM YÜCE \thanks{%
YILDIZ TECHNICAL UNIVERSITY, FACULTY OF ARTS AND SCIENCES, DEPARTMENT OF
MATHEMATICS, 34210, ESENLER, İSTANBUL, TURKEY, ozcanbektas1986@hotmail.com,
sayuce@yildiz.edu.tr.} }
\date{}
\maketitle

\begin{abstract}
In this study, we determine some special Smarandache curves according to
Darboux frame in $E^{3}$ . We give some characterizations and consequences
of Smarandache curves.
\end{abstract}

\textbf{Keywords:} Smarandache curves, Darboux Frame, Normal curvature,
Geodesic \ Curvature, Geodesic torsion.

\textit{2000 AMS Subject Classification:} 53A04

\section{\textbf{Introduction}}

\indent In differential geometry, there are many important consequences and
properties of curves$.~$Researchers follow labours about the curves. In the
light of the existing studies, authors always introduce new curves. Special
Smarandache curves are one of them. This curve is defined as, a regular
curve in Minkowski space-time, whose position vector is composed by Frenet
frame vectors on another regular curve, is called a Smarandache Curve [6].
Special Smarandache curves have been studied by some authors [1, 3, 6]. M.
Turgut and S. Yilmaz have defined a special case of such curves and call it
Smarandache $\mathbf{TB}_{\mathbf{2}}$\emph{\ }Curves in the space $%
E_{1}^{4} $ [6]. They have dealed with a special Smarandache curves which is
defined by the tangent and second binormal vector fields. Besides, they have
computed formulas of this kind curves by the method expressed in [6]. A. T.
Ali has introduced some special Smarandache curves in the Euclidean space
[1]. Special Smarandache curves such as Smarandache curves $\mathbf{TN}_{1}$%
, $\mathbf{TN}_{2}$, $\mathbf{N}_{1}\mathbf{N}_{2}$ and $\mathbf{TN}_{1}%
\mathbf{N}_{2}$ according to Bishop frame in Euclidean 3-space have been
investigated by Çetin et al [3]. Furthermore, they studied differential
geometric properties of these special curves and they calculated first and
second curvature (natural curvatures) of these curves. Also they found the
centers of the curvature spheres and osculating spheres of Smarandache
curves.\newline

\indent In this study, we investigate special Smarandache curves such as
Smarandache $\mathbf{Tg}$, $\mathbf{Tn}$, $\mathbf{gn}$ and $\mathbf{Tgn}$
according to Darboux frame in Euclidean 3-space. Furthermore, we find some
properties of these special curves and we calculate normal curvature,
geodesic curvature and geodesic torsion of these curves.${\,\,\,\,\,\,\,}$

\section{Preliminaries}

\noindent In this section, we give an information about special Smarandache
curves and Darboux frame. Let $M$ ~be an oriented surface and let consider a
curve $\alpha \left( s\right) $ on the surface $M$ . Since the curve $\alpha
\left( s\right) $ is also in space, there exists Frenet frame $\left\{ 
\mathbf{T},\mathbf{N},\mathbf{B}\right\} ~\ $at each points of the curve
where $\mathbf{T~}$is unit tangent vector, $\mathbf{N}~$is principal normal
vector and $\mathbf{B}~$is binormal vector, respectively. The Frenet
equations of the curve $\alpha \left( s\right) $ is given by

\begin{equation*}
\left\{ 
\begin{array}{l}
\mathbf{T}^{^{\prime }}{=\kappa }\mathbf{N} \\ 
\mathbf{N}^{^{\prime }}{=-\kappa \mathbf{T}+\tau }\mathbf{B} \\ 
\mathbf{B}^{^{\prime }}{=-\tau }\mathbf{N}%
\end{array}%
\right.
\end{equation*}%
where ${\kappa }$ and ${\tau }$ are curvature and torsion of the curve $%
\alpha \left( s\right) $ , respectively. Here and in the following, we use
\textquotedblleft dot\textquotedblright\ to denote the derivative with
respect to the arc length parameter of a curve. Since the curve $\alpha
\left( s\right) $ lies on the surface $M$ there exists another frame of the
curve $\alpha \left( s\right) $ which is called Darboux frame and denoted by 
$\left\{ \mathbf{T},\mathbf{g},\mathbf{n}\right\} $. In this frame $\mathbf{T%
}{~}$is the unit tangent of the curve, $\mathbf{n}{~}$is the unit normal of
the surface $M$ \ and $\mathbf{g}{~}$is a unit vector given by $\mathbf{g}{=%
\mathbf{n}\times }\mathbf{T}$. Since the unit tangent $\mathbf{T}{~}$ is
common in both Frenet frame and Darboux frame, the vectors$~\mathbf{N},%
\mathbf{B}{,}\mathbf{g},~\mathbf{n}$ lie on the same plane. So that the
relations between these frames can be given as follows

\begin{equation}
\left[ 
\begin{array}{c}
\mathbf{T} \\ 
\mathbf{g} \\ 
\mathbf{n}%
\end{array}%
\right] =\left[ 
\begin{array}{ccc}
1 & 0 & 0 \\ 
0 & \cos \varphi & \sin \varphi \\ 
0 & -\sin \varphi & \cos \varphi%
\end{array}%
\right] \cdot \left[ 
\begin{array}{c}
{\mathbf{T}} \\ 
\mathbf{N} \\ 
\mathbf{B}%
\end{array}%
\right]  \label{2.1}
\end{equation}%
where$~\varphi $ is the angle between the vectors $\mathbf{g}{~}$and $%
\mathbf{n}$. The derivative formulae of the Darboux frame is

\begin{equation}
\left[ 
\begin{array}{c}
\overset{\cdot }{\mathbf{T}} \\ 
\overset{\cdot }{\mathbf{g}} \\ 
\overset{\cdot }{\mathbf{n}}%
\end{array}%
\right] =\left[ 
\begin{array}{ccc}
0 & k_{g} & k_{n} \\ 
-k_{g} & 0 & \tau _{g} \\ 
-k_{n} & -\tau _{g} & 0%
\end{array}%
\right] \cdot \left[ 
\begin{array}{c}
{\mathbf{T}} \\ 
\mathbf{g} \\ 
\mathbf{n}%
\end{array}%
\right]  \label{2.2}
\end{equation}%
where , $k_{g}~$is the geodesic curvature, $k_{n}~$is the normal curvature
and $\tau _{g}$ is the geodesic torsion of $\alpha \left( s\right) .$

\indent The relations between geodesic curvature, normal curvature, geodesic
torsion and $\kappa $ , $\tau $ are given as follows

\begin{equation}
k_{g}=\kappa \cos \varphi ,~k_{n}=\kappa \sin \varphi \mathrm{\;},\,\mathrm{%
\;}\tau _{g}=\tau +\frac{d\varphi }{ds}.  \label{2.3}
\end{equation}

In the differential geometry of surfaces, for a curve $\alpha \left(
s\right) $ lying on a surface $M$ the followings are well-known

\bigskip i) $\alpha \left( s\right) $ is a geodesic curve $\Leftrightarrow
k_{g}=0,$

ii) $\alpha \left( s\right) $ is an asymptotic line $\Leftrightarrow
k_{n}=0, $

iii) $\alpha \left( s\right) $ is a principal line $\Leftrightarrow \tau
_{g}=0,$ [7].

\noindent \indent Let $\alpha =\alpha \left( s\right) ~$and $\beta =\beta
\left( s^{\ast }\right) ~$be a unit speed regular curves in $E^{3}$ and $%
\left\{ \mathbf{T},\mathbf{N},\mathbf{B}\right\} $ be moving Serret-Frenet
frame of $\alpha \left( s\right) $ . Smarandache curves $\mathbf{TN}$\ are
defined by $\beta \left( s^{\ast }\right) =\frac{1}{\sqrt{2}}\left( \mathbf{T%
}{+}\mathbf{N}\right) $, Smarandache curves $\mathbf{NB}$ are defined by $%
\beta \left( s^{\ast }\right) =\frac{1}{\sqrt{2}}\left( \mathbf{N+B}\right)
~ $and Smarandache curves $\mathbf{TNB~}$are defined by $\beta \left(
s^{\ast }\right) =\frac{1}{\sqrt{3}}\left( \mathbf{T+N+B}\right) .$

\section{Special Smarandache Curves According To Darboux Frame In $E^{3}$}

\indent In this section, we investigate special Smarandache curves according
to Darboux frame in $E^{3}.~$Let $\alpha =\alpha \left( s\right) $ and $%
\beta =\beta \left( s^{\ast }\right) ~$be a unit speed regular curves in $%
E^{3}$ and defined by $\left\{ \mathbf{T},\mathbf{g},\mathbf{n}\right\} ~$%
and $\left\{ \mathbf{T}^{\ast },\mathbf{g}^{\ast },\mathbf{n}^{\ast
}\right\} ~$be Darboux frame of these curves, respectively.

\subsection{$\mathbf{Tg}{-~}$Smarandache Curves}

\textbf{Definition: }Let $M$ ~be an oriented surface in $E^{3}$ and let
consider the arc - length parameter curve $\alpha =\alpha \left( s\right) $
lying fully on $M$. Denote the Darboux frame of $\alpha \left( s\right) $ $%
\left\{ \mathbf{T},\mathbf{g},\mathbf{n}\right\} .$

\indent$~\mathbf{Tg}$- Smarandache curve can be defined by 
\begin{equation}
\beta \left( s^{\ast }\right) =\frac{1}{\sqrt{2}}\left( \mathbf{T+g}\right) .
\label{3.1}
\end{equation}%
\noindent \indent Now, we can investigate Darboux invariants of $\mathbf{Tg}%
-~$Smarandache curve according to $\alpha =\alpha \left( s\right) .$
Differentiating (3.1) with respect to $s$, we get 
\begin{equation}
\beta ^{^{\prime }}=\frac{d\beta }{ds^{\ast }}\frac{ds^{\ast }}{ds}=\frac{-1%
}{\sqrt{2}}\left( k_{g}\mathbf{T}{-k_{g}\mathbf{g}-}\left( k_{n}+\tau {_{g}}%
\right) \mathbf{n}\right) ,  \label{3.2}
\end{equation}%
and 
\begin{equation*}
\mathbf{T}^{\ast }\frac{ds^{\ast }}{ds}=\frac{-1}{\sqrt{2}}\left( k_{g}%
\mathbf{T}{-k_{g}\mathbf{g}-}\left( k_{n}+\tau {_{g}}\right) \mathbf{n}%
\right)
\end{equation*}%
where 
\begin{equation}
\frac{ds^{\ast }}{ds}=\sqrt{\frac{2k_{g}^{2}+\left( k_{n}+\tau {_{g}}\right)
^{2}}{2}}.  \label{3.3}
\end{equation}%
The tangent vector of curve $\beta ~$can be written as follow, 
\begin{equation}
\mathbf{T}^{\ast }=\frac{-1}{\sqrt{2k_{g}^{2}+\left( k_{n}+\tau {_{g}}%
\right) ^{2}}}\left( k_{g}\mathbf{T}{-k_{g}\mathbf{g}-}\left( k_{n}+\tau {%
_{g}}\right) \mathbf{n}\right) .  \label{3.4}
\end{equation}%
Differentiating (3.4) with respect to $s$, we obtain 
\begin{equation}
\frac{d\mathbf{T}^{\ast }}{ds^{\ast }}\frac{ds^{\ast }}{ds}=\frac{-1}{\left(
2k_{g}^{2}+\left( k_{n}+\tau {_{g}}\right) ^{2}\right) ^{\frac{3}{2}}}\left(
\Gamma _{1}\mathbf{T}{+\Gamma _{2}\mathbf{g}+\Gamma _{3}}\mathbf{n}\right)
\label{3.5}
\end{equation}%
where

\bigskip 
\begin{equation*}
\left\{ 
\begin{array}{l}
\Gamma _{1}=\left( k_{n}+\tau {_{g}}\right) \left\{ k_{g}\left(
k_{n}^{^{\prime }}+\tau {_{g}^{^{\prime }}}-k_{n}k_{g}-\tau {_{g}}%
k_{g}\right) -k_{g}^{^{\prime }}\left( k_{n}+\tau {_{g}}\right) -k_{n}\left[
2k_{g}^{2}+\left( k_{n}+\tau {_{g}}\right) ^{2}\right] \right\} -2k_{g}^{4}
\\ 
\Gamma _{2}=\left( k_{n}+\tau {_{g}}\right) \left\{ -k_{g}\left(
k_{n}^{^{\prime }}+\tau {_{g}^{^{\prime }}}+k_{n}k_{g}+\tau {_{g}}%
k_{g}\right) +k_{g}^{^{\prime }}\left( k_{n}+\tau {_{g}}\right) -\tau {_{g}}%
\left[ 2k_{g}^{2}+\left( k_{n}+\tau {_{g}}\right) ^{2}\right] \right\}
-2k_{g}^{4} \\ 
\Gamma _{3}=k_{g}\left( k_{n}+\tau {_{g}}\right) \left[ -2k_{g}^{^{\prime
}}-k_{n}+\tau {_{g}}\left( k_{n}+\tau {_{g}}\right) \right]
+2k_{g}^{2}\left( k_{g}\tau {_{g}+}k_{n}^{^{\prime }}+\tau {_{g}^{^{\prime }}%
}\right) .%
\end{array}%
\right.
\end{equation*}%
Substituting (3.3) in (3.5), we get 
\begin{equation*}
\overset{\cdot }{\mathbf{T}^{\ast }}=\frac{\sqrt{2}}{\left(
2k_{g}^{2}+\left( k_{n}+\tau {_{g}}\right) ^{2}\right) ^{2}}\left( \Gamma
_{1}\mathbf{T}{+\Gamma _{2}\mathbf{g}+\Gamma _{3}}\mathbf{n}\right) .
\end{equation*}%
Then, the curvature and principal normal vector field of curve $\beta ~$are
respectively, 
\begin{equation*}
\kappa ^{\ast }=\left\Vert \overset{\cdot }{\mathbf{T}^{\ast }}\right\Vert =%
\frac{\sqrt{2\left( \Gamma _{1}^{2}{+\Gamma _{2}^{2}+\Gamma _{3}^{2}}\right) 
}}{\left( 2k_{g}^{2}+\left( k_{n}+\tau {_{g}}\right) ^{2}\right) ^{2}}
\end{equation*}%
\noindent and 
\begin{equation*}
\mathbf{N}^{\ast }=\frac{1}{\sqrt{\Gamma _{1}^{2}{+\Gamma _{2}^{2}+\Gamma
_{3}^{2}}}}\left( \Gamma _{1}\mathbf{T}{+\Gamma _{2}\mathbf{g}+\Gamma _{3}}%
\mathbf{n}\right) .
\end{equation*}%
\qquad \qquad \qquad On the other hand, we express 
\begin{equation*}
\mathbf{B}^{\ast }=\mathbf{T}^{\ast }{\times }\mathbf{N}^{\ast }=\frac{1}{%
\sqrt{2k_{g}^{2}+\left( k_{n}+\tau {_{g}}\right) ^{2}}\sqrt{\Gamma _{1}^{2}{%
+\Gamma _{2}^{2}+\Gamma _{3}^{2}}}}\left\vert 
\begin{array}{ccc}
\mathbf{T} & {\mathbf{g}} & \mathbf{n} \\ 
-k_{g} & k_{g} & k_{n}+\tau {_{g}} \\ 
\Gamma _{1} & {\Gamma _{2}} & {\Gamma _{3}}%
\end{array}%
\right\vert .
\end{equation*}%
\qquad So, the binormal vector of curve $\beta ~$is

\begin{equation*}
\mathbf{B}^{\ast }=\frac{1}{\sqrt{2k_{g}^{2}+\left( k_{n}+\tau {_{g}}\right)
^{2}}\sqrt{\Gamma _{1}^{2}{+\Gamma _{2}^{2}+\Gamma _{3}^{2}}}}\left( \mu _{1}%
\mathbf{T}{+\mu _{2}\mathbf{g}+\mu _{3}}\mathbf{n}\right)
\end{equation*}

where 
\begin{equation*}
\left\{ 
\begin{array}{l}
\mu _{1}=k_{g}{\Gamma _{3}-}\left( k_{n}+\tau {_{g}}\right) {\Gamma _{2}} \\ 
\mu _{2}=\left( k_{n}+\tau {_{g}}\right) {\Gamma _{1}+}k_{g}{\Gamma _{3}} \\ 
\mu _{3}=-k_{g}{\Gamma _{2}-}k_{g}{\Gamma _{1.}}%
\end{array}%
\right.
\end{equation*}%
We differentiate (3.2) with respect to s in order to calculate the torsion

\begin{equation*}
\beta ^{^{\prime \prime }}=\frac{-1}{\sqrt{2}}\left\{ \left( k_{g}^{^{\prime
}}+k_{g}^{2}+k_{n}\left( k_{n}+\tau {_{g}}\right) \right) \mathbf{T}{+}%
\left( -k_{g}^{^{\prime }}+k_{g}^{2}+\tau {_{g}}\left( k_{n}+\tau {_{g}}%
\right) \right) {\mathbf{g}+}\left( k_{n}k_{g}-k_{g}\tau {_{g}-}%
k_{g}^{^{\prime }}+\tau {_{g}^{^{\prime }}}\right) \mathbf{n}\right\} ,
\end{equation*}%
and similarly 
\begin{equation*}
\beta ^{^{\prime \prime \prime }}=\frac{-1}{\sqrt{2}}\left( \eta _{1}\mathbf{%
T}{+\eta _{2}\mathbf{g}+\eta _{3}}\mathbf{n}\right) ,
\end{equation*}%
where 
\begin{equation*}
\left\{ 
\begin{array}{l}
\eta _{1}=k_{g}^{^{\prime \prime }}+2k_{n}\left( k_{n}^{^{\prime }}+\tau {%
_{g}^{^{\prime }}}\right) +k_{g}\left( 3k_{g}^{^{\prime }}-\tau {_{g}^{2}}%
-k_{n}^{2}-k_{g}^{2}\right) +k_{n}^{^{\prime }}\left( k_{n}+\tau {_{g}}%
\right) \\ 
\eta _{2}=-k_{g}^{^{\prime \prime }}+2\tau {_{g}}\left( k_{n}^{^{\prime
}}+\tau {_{g}^{^{\prime }}}\right) +k_{g}\left( 3k_{g}^{^{\prime }}+\tau {%
_{g}^{2}}+k_{n}^{2}+k_{g}^{2}\right) +\tau {_{g}^{^{\prime }}}\left(
k_{n}+\tau {_{g}}\right) \\ 
\eta _{3}=-k_{g}^{^{\prime \prime }}-\tau {_{g}^{^{\prime \prime }}}%
+k_{g}\left( k_{n}^{^{\prime }}-\tau {_{g}^{^{\prime }}}\right) +\left(
k_{n}+\tau {_{g}}\right) \left( \tau {_{g}^{2}}+k_{n}^{2}+k_{g}^{2}\right)
+2k_{g}^{^{\prime }}\left( k_{n}-\tau {_{g}}\right) .%
\end{array}%
\right.
\end{equation*}%
The torsion of curve $\beta ~$is

\begin{equation*}
\tau ^{\ast }=\frac{-1}{\sqrt{2}}\frac{\left( \eta _{1}+\eta _{2}\right) %
\left[ k_{g}\left( k_{n}^{^{\prime }}+\tau {_{g}^{^{\prime }}}\right)
-k_{g}^{^{\prime }}\left( k_{n}+\tau {_{g}}\right) \right] +\left(
2k_{g}^{2}+\left( k_{n}+\tau {_{g}}\right) ^{2}\right) \left[ \tau {_{g}}%
\eta _{1}-k_{n}\eta _{2}+k_{g}\eta _{3}\right] }{%
\begin{array}{l}
2\left( k_{g}^{^{\prime }2}+k_{g}^{4}\right) +\left( k_{n}+\tau {_{g}}%
\right) \left[ \left( k_{n}+\tau {_{g}}\right) \left( k_{n}^{2}+\tau {%
_{g}^{2}}\right) +2k_{n}\left( k_{g}^{^{\prime }}+k_{g}^{2}\right) +2\tau {%
_{g}}\left( k_{g}^{^{\prime }}-k_{g}^{2}\right) \right] + \\ 
\left( k_{n}^{^{\prime }}+\tau {_{g}^{^{\prime }}}\right) \left[
k_{n}^{^{\prime }}+\tau {_{g}^{^{\prime }}}-2k_{g}\left( k_{n}-\tau {_{g}}%
\right) \right] +k_{g}^{2}\left( k_{n}-\tau {_{g}}\right) ^{2}%
\end{array}%
}.
\end{equation*}%
The unit normal vector of surface $M$ and unit vector ${\mathbf{g}~}$of $%
\beta $ are as follow. Then, from (2.1) we obtain%
\begin{equation*}
{\mathbf{g}}^{\ast }=\frac{1}{\sqrt{\varsigma }\varepsilon }\left\{ \left( 
\sqrt{\varsigma }\cos \varphi ^{\ast }\Gamma _{1}+\sin \varphi ^{\ast }\mu
_{1}\right) \mathbf{T}{+}\left( \sqrt{\varsigma }\cos \varphi ^{\ast }\Gamma
_{2}+\sin \varphi ^{\ast }\mu _{2}\right) {\mathbf{g}+}\left( \sqrt{%
\varsigma }\cos \varphi ^{\ast }\Gamma _{3}+\sin \varphi ^{\ast }\mu
_{3}\right) \mathbf{n}\right\} ,
\end{equation*}%
and%
\begin{equation*}
\mathbf{n}^{\ast }=\frac{1}{\sqrt{\varsigma }\varepsilon }\left\{ \left( \mu
_{1}\cos \varphi ^{\ast }-\sqrt{\varsigma }\sin \varphi ^{\ast }\Gamma
_{1}\right) \mathbf{T}{+}\left( \mu _{2}\cos \varphi ^{\ast }-\sqrt{%
\varsigma }\sin \varphi ^{\ast }\Gamma _{2}\right) {\mathbf{g}+}\left( \mu
_{3}\cos \varphi ^{\ast }-\sqrt{\varsigma }\sin \varphi ^{\ast }\Gamma
_{3}\right) \mathbf{n}\right\} .
\end{equation*}%
where $\varepsilon =\sqrt{\Gamma _{1}^{2}{+\Gamma _{2}^{2}+\Gamma _{3}^{2}}}%
,~\varsigma =2k_{g}^{2}+\left( k_{n}+\tau {_{g}}\right) ^{2}$ and $\varphi
^{\ast }~$is the angle between the vectors ${\mathbf{g}}^{\ast }~$and $%
\mathbf{n}^{\ast }.$\thinspace Now, we can calculate geodesic curvature,
normal curvature, geodesic torsion of curve , so from (2.3) we get 
\begin{equation*}
k_{g}^{\ast }=\frac{\sqrt{2\left( \Gamma _{1}^{2}{+\Gamma _{2}^{2}+\Gamma
_{3}^{2}}\right) }}{\left( 2k_{g}^{2}+\left( k_{n}+\tau {_{g}}\right)
^{2}\right) ^{2}}\cos \varphi ^{\ast },
\end{equation*}%
and 
\begin{equation*}
k_{n}^{\ast }=\frac{\sqrt{2\left( \Gamma _{1}^{2}{+\Gamma _{2}^{2}+\Gamma
_{3}^{2}}\right) }}{\left( 2k_{g}^{2}+\left( k_{n}+\tau {_{g}}\right)
^{2}\right) ^{2}}\sin \varphi ^{\ast },
\end{equation*}%
and 
\begin{equation*}
\tau _{g}^{\ast }=\frac{-1}{\sqrt{2}}\frac{-1}{\sqrt{2}}\frac{\left( \eta
_{1}+\eta _{2}\right) \left[ k_{g}\left( k_{n}^{^{\prime }}+\tau {%
_{g}^{^{\prime }}}\right) -k_{g}^{^{\prime }}\left( k_{n}+\tau {_{g}}\right) %
\right] +\left( 2k_{g}^{2}+\left( k_{n}+\tau {_{g}}\right) ^{2}\right) \left[
\tau {_{g}}\eta _{1}-k_{n}\eta _{2}+k_{g}\eta _{3}\right] }{%
\begin{array}{l}
2\left( k_{g}^{^{\prime }2}+k_{g}^{4}\right) +\left( k_{n}+\tau {_{g}}%
\right) \left[ \left( k_{n}+\tau {_{g}}\right) \left( k_{n}^{2}+\tau {%
_{g}^{2}}\right) +2k_{n}\left( k_{g}^{^{\prime }}+k_{g}^{2}\right) +2\tau {%
_{g}}\left( k_{g}^{^{\prime }}-k_{g}^{2}\right) \right] + \\ 
\left( k_{n}^{^{\prime }}+\tau {_{g}^{^{\prime }}}\right) \left[
k_{n}^{^{\prime }}+\tau {_{g}^{^{\prime }}}-2k_{g}\left( k_{n}-\tau {_{g}}%
\right) \right] +k_{g}^{2}\left( k_{n}-\tau {_{g}}\right) ^{2}%
\end{array}%
}+\frac{d\varphi ^{\ast }}{ds^{\ast }}.
\end{equation*}

\begin{corollary}
Consider that $\alpha \left( s\right) ~$is a geodesic curve, then the
following equation holds,
\end{corollary}Consider that $\alpha \left( s\right) ~$is a geodesic curve,
then the following equation holds,

$i)~\kappa ^{\ast }=\frac{\sqrt{2~\left( k_{n}^{2}+\tau {_{g}^{2}}\right) }}{%
\left( k_{n}+\tau {_{g}}\right) },$

$ii)~\tau ^{\ast }=\frac{-1}{\sqrt{2}}\frac{\left( k_{n}+\tau {_{g}}\right)
^{2}+\left( k_{n}+\tau {_{g}}\right) \left( k_{n}^{^{\prime }}~\tau {_{g}}%
-k_{n}\tau {_{g}^{^{\prime }}}\right) }{\left( k_{n}+\tau {_{g}}\right)
^{3}\left( k_{n}^{2}+\tau {_{g}^{2}}\right) +\left( k_{n}^{^{\prime }}+\tau {%
_{g}^{^{\prime }}}\right) ^{2}},$

$iii)~k_{g}^{\ast }=\frac{\sqrt{2~\left( k_{n}^{2}+\tau {_{g}^{2}}\right) }}{%
\left( k_{n}+\tau {_{g}}\right) }\cos \varphi ^{\ast },$

$iv)~k_{n}^{\ast }=\frac{\sqrt{2~\left( k_{n}^{2}+\tau {_{g}^{2}}\right) }}{%
\left( k_{n}+\tau {_{g}}\right) }\sin \varphi ^{\ast },$

$v)~\tau _{g}^{\ast }=\frac{-1}{\sqrt{2}}\frac{\left( k_{n}+\tau {_{g}}%
\right) ^{2}+\left( k_{n}+\tau {_{g}}\right) \left( k_{n}^{^{\prime }}~\tau {%
_{g}}-k_{n}\tau {_{g}^{^{\prime }}}\right) }{\left( k_{n}+\tau {_{g}}\right)
^{3}\left( k_{n}^{2}+\tau {_{g}^{2}}\right) +\left( k_{n}^{^{\prime }}+\tau {%
_{g}^{^{\prime }}}\right) ^{2}}+\frac{d\varphi ^{\ast }}{ds^{\ast }}.$

\subsection{$\mathbf{Tn-~}$Smarandache Curves}

\textbf{Definition: }Let $M$ ~be an oriented surface in $E^{3}$ and let
consider the arc - length parameter curve $\alpha =\alpha \left( s\right) $
lying fully on $M$. Denote the Darboux frame of $\alpha \left( s\right) $ $%
\left\{ \mathbf{T},\mathbf{g},\mathbf{n}\right\} .$

$\mathbf{Tn}$-~Smarandache curve can be defined by 
\begin{equation}
\beta \left( s^{\ast }\right) =\frac{1}{\sqrt{2}}\left( \mathbf{T}{+}\mathbf{%
n}\right) .  \label{3.6}
\end{equation}

\indent Now, we can investigate Darboux invariants of $\mathbf{Tn}-~$%
Smarandache curve according to $\alpha =\alpha \left( s\right) .$
Differentiating (3.6) with respect to $s$, we get 
\begin{equation}
\beta ^{^{\prime }}=\frac{d\beta }{ds^{\ast }}\frac{ds^{\ast }}{ds}=\frac{-1%
}{\sqrt{2}}\left( k_{n}\mathbf{T}{+\left( \tau {_{g}-k_{g}}\right) \mathbf{g}%
-k_{n}}\mathbf{n}\right) ,  \label{3.7}
\end{equation}%
and 
\begin{equation*}
\mathbf{T}^{\ast }\frac{ds^{\ast }}{ds}=\frac{-1}{\sqrt{2}}\left( k_{n}%
\mathbf{T}{+\left( \tau {_{g}-k_{g}}\right) \mathbf{g}-k_{n}}\mathbf{n}%
\right)
\end{equation*}%
where 
\begin{equation}
\frac{ds^{\ast }}{ds}=\sqrt{\frac{2k_{n}^{2}+{\left( \tau {_{g}-k_{g}}%
\right) }^{2}}{2}}.  \label{3.8}
\end{equation}%
The tangent vector of curve $\beta ~$can be written as follow, 
\begin{equation}
\mathbf{T}^{\ast }=\frac{-1}{\sqrt{2k_{n}^{2}+{\left( \tau {_{g}-k_{g}}%
\right) }^{2}}}\left( k_{n}\mathbf{T}{+\left( \tau {_{g}-k_{g}}\right) 
\mathbf{g}-k_{n}}\mathbf{n}\right) .  \label{3.9}
\end{equation}%
Differentiating (3.9) with respect to $s$, we obtain 
\begin{equation}
\frac{d\mathbf{T}^{\ast }}{ds^{\ast }}\frac{ds^{\ast }}{ds}=\frac{-1}{\left(
2k_{n}^{2}+{\left( \tau {_{g}-k_{g}}\right) }^{2}\right) ^{\frac{3}{2}}}%
\left( \gamma _{1}\mathbf{T}{+\gamma _{2}\mathbf{g}+\gamma _{3}}\mathbf{n}%
\right)  \label{3.10}
\end{equation}%
where 
\begin{equation}
\left\{ 
\begin{array}{l}
\gamma _{1}={\left( \tau {_{g}-k_{g}}\right) }\left\{ k_{n}\left(
-k_{g}^{^{\prime }}+\tau {_{g}^{^{\prime }}}+k_{n}k_{g}-\tau {_{g}}%
k_{n}\right) -k_{n}^{^{\prime }}{\left( \tau {_{g}-k_{g}}\right) }+k_{g}%
\left[ 2k_{n}^{2}+\tau {_{g}-k_{g}}\right] \right\} -2k_{n}^{4} \\ 
\gamma _{2}=k_{n}{\left( \tau {_{g}-k_{g}}\right) }\left[ 2k_{n}^{^{\prime
}}+\tau {_{g}}\left( k_{g}^{2}-\tau {_{g}^{2}}\right) \right]
-2k_{n}^{2}\left( k_{n}k_{g}+\tau {_{g}^{^{\prime }}-}k_{g}^{^{\prime
}}-k_{n}\tau {_{g}}\right) \\ 
\gamma _{3}={\left( \tau {_{g}-k_{g}}\right) }\left\{ -k_{n}\left(
-k_{g}^{^{\prime }}+\tau {_{g}^{^{\prime }}}-k_{n}k_{g}+\tau {_{g}}%
k_{n}\right) +k_{n}^{^{\prime }}\left( \tau {_{g}-k_{g}}\right) -\tau {_{g}}%
\left[ 2k_{n}^{2}+\tau {_{g}-k_{g}}\right] \right\} -2k_{n}^{4}.%
\end{array}%
\right.  \notag
\end{equation}%
Substituting (3.8) in (3.10), we get 
\begin{equation*}
\overset{\cdot }{\mathbf{T}^{\ast }}=\frac{\sqrt{2}}{\left( 2k_{n}^{2}+{%
\left( \tau {_{g}-k_{g}}\right) }^{2}\right) ^{2}}\left( \gamma _{1}\mathbf{T%
}{+\gamma _{2}\mathbf{g}+\gamma _{3}}\mathbf{n}\right) .
\end{equation*}%
Then, the curvature and principal normal vector field of curve $\beta ~$are
respectively, 
\begin{equation*}
\kappa ^{\ast }=\left\Vert \overset{\cdot }{\mathbf{T}^{\ast }}\right\Vert =%
\frac{\sqrt{2~\left( \gamma _{1}^{2}{+\gamma _{2}^{2}+\gamma _{3}^{2}}%
\right) }}{\left( 2k_{n}^{2}+{\left( \tau {_{g}-k_{g}}\right) }^{2}\right)
^{2}}
\end{equation*}%
\noindent and 
\begin{equation*}
\mathbf{N}^{\ast }=\frac{1}{\sqrt{\gamma _{1}^{2}{+\gamma _{2}^{2}+\gamma
_{3}^{2}}}}\left( \gamma _{1}\mathbf{T}{+\gamma _{2}\mathbf{g}+\gamma _{3}}%
\mathbf{n}\right) .
\end{equation*}%
\qquad \qquad \qquad On the other hand, we express 
\begin{equation*}
\mathbf{B}^{\ast }=\mathbf{T}^{\ast }{\times }\mathbf{N}^{\ast }=\frac{1}{%
\sqrt{2k_{n}^{2}+{\left( \tau {_{g}-k_{g}}\right) }^{2}}\sqrt{\gamma _{1}^{2}%
{+\gamma _{2}^{2}+\gamma _{3}^{2}}}}\left\vert 
\begin{array}{ccc}
\mathbf{T} & {\mathbf{g}} & \mathbf{n} \\ 
-k_{n} & {k_{g}-}\tau {_{g}} & k_{n} \\ 
\gamma _{1} & \gamma {_{2}} & \gamma {_{3}}%
\end{array}%
\right\vert .
\end{equation*}%
\qquad So, the binormal vector of curve $\beta ~$is

\begin{equation*}
\mathbf{B}^{\ast }=\frac{1}{\sqrt{2k_{n}^{2}+{\left( \tau {_{g}-k_{g}}%
\right) }^{2}}\sqrt{\gamma _{1}^{2}{+\gamma _{2}^{2}+\gamma _{3}^{2}}}}%
\left( \nu _{1}\mathbf{T}{+\nu _{2}\mathbf{g}+\nu _{3}}\mathbf{n}\right)
\end{equation*}

where 
\begin{equation*}
\left\{ 
\begin{array}{l}
\nu _{1}=\left( {k_{g}-}\tau {_{g}}\right) {\gamma _{3}-k_{n}\gamma _{2}} \\ 
\nu _{2}={k_{n}\gamma _{1}+}k_{n}{\gamma _{3}} \\ 
\nu _{3}=-k_{n}{\gamma _{2}+\left( \tau {_{g}-k_{g}}\right) }\gamma _{1}.%
\end{array}%
\right.
\end{equation*}%
We differentiate (3.7) with respect to $s$ in order to calculate the torsion

\begin{equation*}
\beta ^{^{\prime \prime }}=\frac{-1}{\sqrt{2}}\left\{ \left( k_{g}^{^{\prime
}}+k_{n}^{2}-k_{g}\left( \tau {_{g}-k_{g}}\right) \right) \mathbf{T}{+}%
\left( k_{n}k_{g}+\left( \tau {_{g}^{^{\prime }}-}k_{g}^{^{\prime }}\right)
+k_{n}\tau {_{g}}\right) {\mathbf{g}+}\left( k_{n}^{2}+\left( \tau {%
_{g}-k_{g}}\right) \tau {_{g}}-k_{n}^{^{\prime }}\right) \mathbf{n}\right\} ,
\end{equation*}%
and similarly 
\begin{equation*}
\beta ^{^{\prime \prime \prime }}=\frac{-1}{\sqrt{2}}\left( \omega _{1}%
\mathbf{T}{+\omega _{2}\mathbf{g}+\omega _{3}}\mathbf{n}\right) ,
\end{equation*}%
where 
\begin{equation*}
\left\{ 
\begin{array}{l}
\omega _{1}=k_{n}^{^{\prime \prime }}-2k_{g}\left( \tau {_{g}^{^{\prime }}-}%
k_{g}^{^{\prime }}\right) +k_{n}\left( 3k_{n}^{^{\prime }}-\tau {_{g}^{2}}%
-k_{n}^{2}-k_{g}^{2}\right) +k_{g}^{^{\prime }}\left( k_{g}-\tau {_{g}}%
\right) \\ 
\omega _{2}=k_{g}^{^{\prime \prime }}-\tau {_{g}^{^{\prime \prime }}}%
+k_{n}\left( k_{g}^{^{\prime }}+\tau {_{g}^{^{\prime }}}\right) +\left(
k_{g}-\tau {_{g}}\right) \left( \tau {_{g}^{2}}+k_{n}^{2}+k_{g}^{2}\right)
+\tau {_{g}}\left( 2k_{n}^{^{\prime }}-k_{n}^{2}\right)
+2k_{g}k_{n}^{^{\prime }} \\ 
\omega _{3}=-k_{n}^{^{\prime \prime }}+2\tau {_{g}}\left( \tau {%
_{g}^{^{\prime }}-}k_{g}^{^{\prime }}\right) +k_{n}\left( 3k_{n}^{^{\prime
}}+\tau {_{g}^{2}}+\tau {_{g}}\left( k_{g}+\tau {_{g}}\right) \right)
+\left( k_{g}-\tau {_{g}}\right) \left( k_{n}k_{g}-\tau {_{g}^{^{\prime }}}%
\right) .%
\end{array}%
\right.
\end{equation*}%
The torsion of curve $\beta ~$is

\begin{equation*}
\tau ^{\ast }=\frac{-1}{\sqrt{2}}\frac{%
\begin{array}{l}
\left( \omega _{1}+\omega _{3}\right) \left[ k_{n}\left( \tau {%
_{g}^{^{\prime }}}-k_{g}^{^{\prime }}\right) -k_{n}^{^{\prime }}\left( \tau {%
_{g}-}k_{g}\right) +k_{n}^{2}\left( \tau {_{g}+}k_{g}\right) \right] +\left(
\tau {_{g}-}k_{g}\right) ^{2}\left( \tau {_{g}}\omega _{1}+k_{g}\omega
_{3}\right) + \\ 
k_{n}\left( \tau {_{g}-}k_{g}\right) \left[ k_{n}\left( \omega _{1}-\omega
_{3}\right) +\omega _{2}\left( k_{g}-\tau {_{g}}\right) \right]
-2k_{n}\left( k_{n}^{^{\prime }}+k_{n}^{2}\right) \omega _{2}%
\end{array}%
}{%
\begin{array}{l}
2\left( k_{n}^{^{\prime }2}+k_{n}^{4}\right) +\left( \tau {_{g}-}%
k_{g}\right) \left[ 2\left( k_{n}^{2}-k_{n}^{^{\prime }}\right) -2\left(
k_{n}^{2}-k_{n}^{^{\prime }}\right) k_{g}\right] +\left( \tau {_{g}-}%
k_{g}\right) ^{2}\left( 1+k_{g}^{2}\right) + \\ 
\left( \tau {_{g}}-k_{g}^{^{\prime }}\right) \left[ \left( \tau {%
_{g}^{^{\prime }}}-k_{g}^{^{\prime }}\right) +2\left( k_{n}k_{g}+k_{n}\tau {%
_{g}}\right) \right] +k_{n}^{2}\left( k_{g}+\tau {_{g}}\right) ^{2}%
\end{array}%
}.
\end{equation*}%
The unit normal vector of surface $M$ and unit vector $\mathbf{g~}$of $\beta 
$ are as follow. Then, from (2.1) we obtain%
\begin{equation*}
{\overrightarrow{g}}^{\ast }=\frac{1}{\sqrt{\psi }\pi }\left\{ \left( \sqrt{%
\psi }\cos \varphi ^{\ast }\gamma _{1}+\sin \varphi ^{\ast }\nu _{1}\right) 
\mathbf{T}{+}\left( \sqrt{\psi }\cos \varphi ^{\ast }\gamma _{2}+\sin
\varphi ^{\ast }\nu _{2}\right) {\mathbf{g}+}\left( \sqrt{\psi }\cos \varphi
^{\ast }\gamma _{3}+\sin \varphi ^{\ast }\nu _{3}\right) \mathbf{n}\right\} ,
\end{equation*}%
and%
\begin{equation*}
{\overrightarrow{N}}^{\ast }=\frac{1}{\sqrt{\psi }\pi }\left\{ \left( \nu
_{1}\cos \varphi ^{\ast }-\sqrt{\psi }\sin \varphi ^{\ast }\gamma
_{1}\right) {\overrightarrow{t}+}\left( \nu _{2}\cos \varphi ^{\ast }-\sqrt{%
\psi }\sin \varphi ^{\ast }\gamma _{2}\right) {\overrightarrow{g}+}\left(
\nu _{3}\cos \varphi ^{\ast }-\sqrt{\psi }\sin \varphi ^{\ast }\gamma
_{3}\right) {\overrightarrow{N}}\right\} .
\end{equation*}%
where $\pi =\sqrt{\gamma _{1}^{2}{+\gamma _{2}^{2}+\gamma _{3}^{2}}},~\psi
=2k_{n}^{2}+{\left( \tau {_{g}-k_{g}}\right) }^{2}$ and $\varphi ^{\ast }~$%
is the angle between the vectors ${\mathbf{g}}^{\ast }~$and $\mathbf{n}%
^{\ast }.$Now, we can calculate geodesic curvature, normal curvature,
geodesic torsion of curve , so from (2.3) we get 
\begin{equation*}
k_{g}^{\ast }=\frac{\sqrt{2\left( \gamma _{1}^{2}{+\gamma _{2}^{2}+\gamma
_{3}^{2}}\right) }}{\left( 2k_{n}^{2}+\left( \tau {_{g}-}k_{g}\right)
^{2}\right) ^{2}}\cos \varphi ^{\ast },
\end{equation*}%
and 
\begin{equation*}
k_{n}^{\ast }=\frac{\sqrt{2\left( \gamma _{1}^{2}{+\gamma _{2}^{2}+\gamma
_{3}^{2}}\right) }}{\left( 2k_{n}^{2}+\left( \tau {_{g}-}k_{g}\right)
^{2}\right) ^{2}}\sin \varphi ^{\ast },
\end{equation*}%
and 
\begin{equation*}
\tau _{g}^{\ast }=\frac{-1}{\sqrt{2}}\frac{%
\begin{array}{l}
\left( \omega _{1}+\omega _{3}\right) \left[ k_{n}\left( \tau {%
_{g}^{^{\prime }}}-k_{g}^{^{\prime }}\right) -k_{n}^{^{\prime }}\left( \tau {%
_{g}-}k_{g}\right) +k_{n}^{2}\left( \tau {_{g}+}k_{g}\right) \right] +\left(
\tau {_{g}-}k_{g}\right) ^{2}\left( \tau {_{g}}\omega _{1}+k_{g}\omega
_{3}\right) + \\ 
k_{n}\left( \tau {_{g}-}k_{g}\right) \left[ k_{n}\left( \omega _{1}-\omega
_{3}\right) +\omega _{2}\left( k_{g}-\tau {_{g}}\right) \right]
-2k_{n}\left( k_{n}^{^{\prime }}+k_{n}^{2}\right) \omega _{2}%
\end{array}%
}{%
\begin{array}{l}
2\left( k_{n}^{^{\prime }2}+k_{n}^{4}\right) +\left( \tau {_{g}-}%
k_{g}\right) \left[ 2\left( k_{n}^{2}-k_{n}^{^{\prime }}\right) -2\left(
k_{n}^{2}-k_{n}^{^{\prime }}\right) k_{g}\right] +\left( \tau {_{g}-}%
k_{g}\right) ^{2}\left( 1+k_{g}^{2}\right) + \\ 
\left( \tau {_{g}}-k_{g}^{^{\prime }}\right) \left[ \left( \tau {%
_{g}^{^{\prime }}}-k_{g}^{^{\prime }}\right) +2\left( k_{n}k_{g}+k_{n}\tau {%
_{g}}\right) \right] +k_{n}^{2}\left( k_{g}+\tau {_{g}}\right) ^{2}%
\end{array}%
}+\frac{d\varphi ^{\ast }}{ds^{\ast }}.
\end{equation*}%
\begin{corollary}
Consider that $\alpha \left( s\right) ~$is a geodesic curve, then the
following equation holds,
\end{corollary}Consider that $\alpha \left( s\right) ~$is an asymptotic
line, then the following equation holds,

$i)~\kappa ^{\ast }=\frac{\sqrt{2~\left( k_{g}^{2}+\tau {_{g}^{2}}\right) }}{%
\left( \tau {_{g}-}k_{g}\right) ^{2}},$

$ii)~\tau ^{\ast }=\frac{-1}{\sqrt{2}}\frac{\left( k_{g}-\tau {_{g}}\right)
\left( k_{g}^{^{\prime }}\tau {_{g}}-k_{g}\tau {_{g}^{^{\prime }}}~~\right) 
}{\left( 1+k_{g}^{2}\right) \left( \tau {_{g}-}k_{g}^{^{\prime }}\right)
\left( \tau {_{g}^{^{\prime }}}-k_{g}^{^{\prime }}\right) \left( \tau {_{g}-}%
k_{g}\right) ^{-2}},$

$iii)~k_{g}^{\ast }=\frac{\sqrt{2~\left( k_{g}^{2}+\tau {_{g}^{2}}\right) }}{%
\left( \tau {_{g}-}k_{g}\right) ^{2}}\cos \varphi ^{\ast },$

$iv)~k_{n}^{\ast }=\frac{\sqrt{2~\left( k_{g}^{2}+\tau {_{g}^{2}}\right) }}{%
\left( \tau {_{g}-}k_{g}\right) ^{2}}\sin \varphi ^{\ast },$

$v)~\tau _{g}^{\ast }=\frac{-1}{\sqrt{2}}\frac{\left( k_{g}-\tau {_{g}}%
\right) \left( k_{g}^{^{\prime }}\tau {_{g}}-k_{g}\tau {_{g}^{^{\prime }}}%
~~~~\right) }{\left( 1+k_{g}^{2}\right) \left( \tau {_{g}-}k_{g}^{^{\prime
}}\right) \left( \tau {_{g}^{^{\prime }}}-k_{g}^{^{\prime }}\right) \left(
\tau {_{g}-}k_{g}\right) ^{-2}}+\frac{d\varphi ^{\ast }}{ds^{\ast }}.$

\subsection{$\mathbf{gn}{-~}$Smarandache Curves}

\textbf{Definition:}Let $M$ ~be an oriented surface in $E^{3}$ and let
consider the arc - length parameter curve $\alpha =\alpha \left( s\right) $
lying fully on $M$. Denote the Darboux frame of $\alpha \left( s\right) $ $%
\left\{ \mathbf{T},\mathbf{g},\mathbf{n}\right\} .$

$\mathbf{gn}$\textbf{\ }- Smarandache curve can be defined by 
\begin{equation}
\beta \left( s^{\ast }\right) =\frac{1}{\sqrt{2}}\left( \mathbf{g}+\mathbf{n}%
\right) .  \label{3.11}
\end{equation}

\indent Now, we can investigate Darboux invariants of $\mathbf{gn}-~$%
Smarandache curve according to $\alpha =\alpha \left( s\right) .$
Differentiating (3.11) with respect to $s$, we get 
\begin{equation}
\beta ^{^{\prime }}=\frac{d\beta }{ds^{\ast }}\frac{ds^{\ast }}{ds}=\frac{-1%
}{\sqrt{2}}\left( \left( k_{n}+{k_{g}}\right) \mathbf{T}{+\tau {_{g}}\mathbf{%
g}-\tau {_{g}}}\mathbf{n}\right) ,  \label{3.12}
\end{equation}%
and 
\begin{equation*}
\mathbf{T}^{\ast }\frac{ds^{\ast }}{ds}=\frac{-1}{\sqrt{2}}\left( \left(
k_{n}+{k_{g}}\right) \mathbf{T}{+\tau {_{g}}\mathbf{g}-\tau {_{g}}}\mathbf{n}%
\right)
\end{equation*}%
where 
\begin{equation}
\frac{ds^{\ast }}{ds}=\sqrt{\frac{2{\tau {_{g}^{2}}}+\left( k_{n}+{k_{g}}%
\right) ^{2}}{2}}.  \label{3.13}
\end{equation}%
The tangent vector of curve $\beta ~$can be written as follow, 
\begin{equation}
\mathbf{T}^{\ast }=\frac{-1}{\sqrt{2{\tau {_{g}^{2}}}+\left( k_{n}+{k_{g}}%
\right) ^{2}}}\left( \left( k_{n}+{k_{g}}\right) \mathbf{T}{+\tau {_{g}}%
\mathbf{g}-\tau {_{g}}}\mathbf{n}\right) .  \label{3.14}
\end{equation}%
Differentiating (3.14) with respect to $s$, we obtain 
\begin{equation}
\frac{d\mathbf{T}^{\ast }}{ds^{\ast }}\frac{ds^{\ast }}{ds}=\frac{-1}{\left(
2{\tau {_{g}^{2}}}+\left( k_{n}+{k_{g}}\right) ^{2}\right) ^{\frac{3}{2}}}%
\left( \lambda _{1}\mathbf{T}{+\lambda _{2}\mathbf{g}+\lambda _{3}}\mathbf{n}%
\right)  \label{3.15}
\end{equation}%
where 
\begin{equation}
\left\{ 
\begin{array}{l}
\lambda _{1}=2\tau {_{g}\tau {_{g}^{^{\prime }}}\left( k_{n}+{k_{g}}\right) -%
}2\tau {_{g}^{2}}\left( k_{n}^{^{\prime }}+k_{g}^{^{\prime }}\right) +\tau {%
_{g}}\left( k_{g}-k_{n}\right) \left( 2\tau {_{g}^{2}+\left( k_{n}+{k_{g}}%
\right) }^{2}\right) \\ 
\lambda _{2}=-2\tau {_{g}^{4}}+\tau {_{g}\left( k_{n}+{k_{g}}\right) }\left[
\left( k_{n}^{^{\prime }}+k_{g}^{^{\prime }}\right) -2\tau {_{g}}k_{g}\right]
-{\left( k_{n}+{k_{g}}\right) }^{2}\left( \left( k_{n}+k_{g}\right) {k_{g}}+{%
\tau {_{g}^{^{\prime }}}+}\tau {_{g}^{2}}\right) \\ 
\lambda _{3}=-2\tau {_{g}^{4}+}\tau {_{g}\left( k_{n}+{k_{g}}\right) }\left[
\left( k_{n}^{^{\prime }}+k_{g}^{^{\prime }}\right) -2\tau {_{g}}k_{n}\right]
+{\left( k_{n}+{k_{g}}\right) }^{2}\left( -\left( k_{n}+k_{g}\right) {k_{n}}+%
{\tau {_{g}^{^{\prime }}}-}\tau {_{g}^{2}}\right)%
\end{array}%
\right.  \notag
\end{equation}%
Substituting (3.13) in (3.15), we get 
\begin{equation*}
\overset{\cdot }{\mathbf{T}^{\ast }}=\frac{\sqrt{2}}{\left( 2\tau {_{g}^{2}}+%
{\left( k_{n}+{k_{g}}\right) }^{2}\right) ^{2}}\left( \lambda _{1}\mathbf{T}{%
+\lambda _{2}\mathbf{g}+\lambda _{3}}\mathbf{n}\right) .
\end{equation*}%
Then, the curvature and principal normal vector field of curve $\beta ~$are
respectively, 
\begin{equation*}
\kappa ^{\ast }=\left\Vert \overset{\cdot }{\mathbf{T}^{\ast }}\right\Vert =%
\frac{\sqrt{2~\left( \lambda _{1}^{2}{+\lambda _{2}^{2}+\lambda _{3}^{2}}%
\right) }}{\left( 2\tau {_{g}^{2}}+{\left( k_{n}+{k_{g}}\right) }^{2}\right)
^{2}}
\end{equation*}%
\noindent and 
\begin{equation*}
\mathbf{N}^{\ast }=\frac{1}{\sqrt{\lambda _{1}^{2}{+\lambda _{2}^{2}+\lambda
_{3}^{2}}}}\left( \lambda _{1}\mathbf{T}{+\lambda _{2}\mathbf{g}+\lambda _{3}%
}\mathbf{n}\right) .
\end{equation*}%
\qquad \qquad \qquad On the other hand, we express 
\begin{equation*}
\mathbf{B}^{\ast }=\mathbf{T}^{\ast }{\times }\mathbf{N}^{\ast }=\frac{1}{%
\sqrt{2k_{n}^{2}+{\left( \tau {_{g}-k_{g}}\right) }^{2}}\sqrt{\lambda
_{1}^{2}{+\lambda _{2}^{2}+\lambda _{3}^{2}}}}\left\vert 
\begin{array}{ccc}
\mathbf{T} & {\mathbf{g}} & \mathbf{n} \\ 
-{\left( k_{n}+{k_{g}}\right) } & {-}\tau {_{g}} & \tau {_{g}} \\ 
\lambda _{1} & \lambda {_{2}} & \lambda {_{3}}%
\end{array}%
\right\vert .
\end{equation*}%
\qquad So, the binormal vector of curve $\beta ~$is

\begin{equation*}
\mathbf{B}^{\ast }=\frac{1}{\sqrt{2\tau {_{g}^{2}}+{\left( k_{n}+{k_{g}}%
\right) }^{2}}\sqrt{\lambda _{1}^{2}{+\lambda _{2}^{2}+\lambda _{3}^{2}}}}%
\left( \rho _{1}\mathbf{T}{+\rho _{2}\mathbf{g}+\rho _{3}}\mathbf{n}\right)
\end{equation*}

where 
\begin{equation*}
\left\{ 
\begin{array}{l}
\rho _{1}={-}\tau {_{g}\lambda _{3}-\tau {_{g}}\lambda _{2}} \\ 
{\rho _{2}}={\tau {_{g}}\lambda _{1}+\left( k_{n}+{k_{g}}\right) \lambda _{3}%
} \\ 
{\rho _{3}}=-{\left( k_{n}+{k_{g}}\right) \lambda _{2}+}\tau {_{g}\lambda
_{1}}.%
\end{array}%
\right.
\end{equation*}%
We differentiate (3.7) with respect to s in order to calculate the torsion

\begin{equation*}
\beta ^{^{\prime \prime }}=\frac{-1}{\sqrt{2}}\left\{ \left( k_{g}^{^{\prime
}}+k_{n}^{^{\prime }}+\tau {_{g}}\left( k_{n}{-k_{g}}\right) \right) \mathbf{%
T}{+}\left( k_{g}{\left( k_{n}+{k_{g}}\right) }+{\tau {_{g}^{^{\prime }}}+}%
\tau {_{g}^{2}}\right) {\mathbf{g}+}\left( k_{n}{\left( k_{n}+{k_{g}}\right) 
}-{\tau {_{g}^{^{\prime }}}+}\tau {_{g}^{2}}\right) \mathbf{n}\right\} ,
\end{equation*}%
and similarly 
\begin{equation*}
\beta ^{^{\prime \prime \prime }}=\frac{-1}{\sqrt{2}}\left( \chi _{1}\mathbf{%
T}{+\chi _{2}\mathbf{g}+\chi _{3}}\mathbf{n}\right) ,
\end{equation*}%
where 
\begin{equation*}
\left\{ 
\begin{array}{l}
\chi _{1}=k_{n}^{^{\prime \prime }}+k_{g}^{^{\prime \prime }}-2\overset{%
\cdot }{\tau {_{g}}}\left( k_{n}-k{_{g}}\right) +\tau {_{g}}\left(
k_{n}^{^{\prime }}-k_{g}^{^{\prime }}\right) -\left( k_{n}+k{_{g}}\right)
\left( \tau {_{g}^{2}}+k_{n}^{2}+k_{g}^{2}\right) \\ 
\chi _{2}={\tau {_{g}^{^{\prime \prime }}}}+2k_{g}\left( k_{n}^{^{\prime
}}+k_{g}^{^{\prime }}\right) +3\tau {_{g}\tau {_{g}^{^{\prime }}}}+\left(
k_{n}+k{_{g}}\right) \left( k_{g}^{^{\prime }}-\tau {_{g}}k_{n}\right)
+k_{g}\tau {_{g}}\left( k_{n}-k{_{g}}\right) -\tau {_{g}^{3}} \\ 
\chi _{3}={\tau {_{g}^{^{\prime \prime }}}}+2k_{n}\left( k_{n}^{^{\prime
}}+k_{g}^{^{\prime }}\right) +3\tau {_{g}\tau {_{g}^{^{\prime }}}}+\left(
k_{n}+k{_{g}}\right) \left( k_{g}^{^{\prime }}+\tau {_{g}}k_{g}\right)
+k_{n}\tau {_{g}}\left( k_{n}-k{_{g}}\right) +\tau {_{g}^{3}}.%
\end{array}%
\right.
\end{equation*}%
The torsion of curve $\beta ~$is $\ $

\bigskip $\tau ^{\ast }=\frac{-1}{\sqrt{2}}\frac{%
\begin{array}{l}
\left( k_{n}+k{_{g}}\right) ^{2}[\left( \left( {\tau {_{g}^{^{\prime }}}}%
+\tau {_{g}}\right) ^{2}+k{_{g}}\left( k_{n}+k{_{g}}\right) \right) \chi
_{3}-\left( \left( {\tau {_{g}^{^{\prime }}}}-\tau {_{g}}\right) ^{2}-k{_{n}}%
\left( k_{n}+k{_{g}}\right) \right) \chi _{2}+ \\ 
\tau {_{g}}\left( k_{n}+k{_{g}}\right) \chi _{1}]-\tau {_{g}}\left( \chi
_{2}+\chi _{3}\right) \left[ \left( k_{n}^{^{\prime }}+k_{g}^{^{\prime
}}\right) +\tau {_{g}}\left( k_{n}-k{_{g}}\right) \right] ^{2}%
\end{array}%
}{%
\begin{array}{l}
\left( k_{n}+k{_{g}}\right) ^{2}\left( k_{n}^{2}+k_{g}^{2}\right) +2\left(
k_{n}+k{_{g}}\right) \left[ k{_{g}}\left( {\tau {_{g}^{^{\prime }}}}+\tau {%
_{g}^{2}}\right) +2k{_{n}}\left( \tau {_{g}^{2}-\tau {_{g}^{^{\prime }}}}%
\right) \right] + \\ 
2{\tau {_{g}^{^{\prime }}}}^{2}+2\tau {_{g}^{4}+}\left[ \left(
k_{n}^{^{\prime }}+k_{g}^{^{\prime }}\right) +\tau {_{g}}\left( k_{n}-k{_{g}}%
\right) \right] ^{2}%
\end{array}%
}.$

\bigskip The unit normal vector of surface $M$ and unit vector $\mathbf{g}{~}
$of $\beta $ are as follow. Then, from (2.1) we obtain%
\begin{equation*}
\mathbf{g}^{\ast }=\frac{1}{\sqrt{\Delta }\Omega }\left\{ \left( \sqrt{%
\Delta }\cos \varphi ^{\ast }{\lambda }_{1}+\sin \varphi ^{\ast }\rho
_{1}\right) \mathbf{T}{+}\left( \sqrt{\Delta }\cos \varphi ^{\ast }{\lambda }%
_{2}+\sin \varphi ^{\ast }\rho _{2}\right) \mathbf{g}{+}\left( \sqrt{\Delta }%
\cos \varphi ^{\ast }\lambda _{3}+\sin \varphi ^{\ast }\rho _{3}\right) 
\mathbf{n}\right\} ,
\end{equation*}%
and%
\begin{equation*}
\mathbf{n}^{\ast }=\frac{1}{\sqrt{\Delta }\Omega }\left\{ \left( \rho
_{1}\cos \varphi ^{\ast }-\sqrt{\Delta }\sin \varphi ^{\ast }{\lambda }%
_{1}\right) \mathbf{T}{+}\left( \rho _{2}\cos \varphi ^{\ast }-\sqrt{\Delta }%
\sin \varphi ^{\ast }{\lambda }_{2}\right) \mathbf{g}{+}\left( \rho _{3}\cos
\varphi ^{\ast }-\sqrt{\Delta }\sin \varphi ^{\ast }\lambda _{3}\right) 
\mathbf{n}\right\} .
\end{equation*}%
where $\Omega =\sqrt{\lambda _{1}^{2}{+\lambda _{2}^{2}+\lambda _{3}^{2}}}%
,~\Delta =2\tau {_{g}^{2}}+{\left( k_{n}+{k_{g}}\right) }^{2}$ and $\varphi
^{\ast }~$is the angle between the vectors ${\mathbf{g}}^{\ast }~$and $%
\mathbf{n}^{\ast }.~$Now, we can calculate geodesic curvature, normal
curvature, geodesic torsion of curve , so from (2.3) we get 
\begin{equation*}
k_{g}^{\ast }=\frac{\sqrt{2\left( \lambda _{1}^{2}{+\lambda _{2}^{2}+\lambda
_{3}^{2}}\right) }}{\left( 2\tau {_{g}^{2}}+{\left( k_{n}+{k_{g}}\right) }%
^{2}\right) ^{2}}\cos \varphi ^{\ast },
\end{equation*}%
and 
\begin{equation*}
k_{n}^{\ast }=\frac{\sqrt{2\left( \lambda _{1}^{2}{+\lambda _{2}^{2}+\lambda
_{3}^{2}}\right) }}{\left( 2\tau {_{g}^{2}}+{\left( k_{n}+{k_{g}}\right) }%
^{2}\right) ^{2}}\sin \varphi ^{\ast },
\end{equation*}%
and

\begin{equation*}
\tau ^{\ast }=\frac{-1}{\sqrt{2}}\frac{%
\begin{array}{l}
\left( k_{n}+k{_{g}}\right) ^{2}[\left( \left( {\tau {_{g}^{^{\prime }}}}%
+\tau {_{g}}\right) ^{2}+k{_{g}}\left( k_{n}+k{_{g}}\right) \right) \chi
_{3}-\left( \left( {\tau {_{g}^{^{\prime }}}}-\tau {_{g}}\right) ^{2}-k{_{n}}%
\left( k_{n}+k{_{g}}\right) \right) \chi _{2}+ \\ 
\tau {_{g}}\left( k_{n}+k{_{g}}\right) \chi _{1}]-\tau {_{g}}\left( \chi
_{2}+\chi _{3}\right) \left[ \left( k_{n}^{^{\prime }}+k_{g}^{^{\prime
}}\right) +\tau {_{g}}\left( k_{n}-k{_{g}}\right) \right] ^{2}%
\end{array}%
}{%
\begin{array}{l}
\left( k_{n}+k{_{g}}\right) ^{2}\left( k_{n}^{2}+k_{g}^{2}\right) +2\left(
k_{n}+k{_{g}}\right) \left[ k{_{g}}\left( {\tau {_{g}^{^{\prime }}}}+\tau {%
_{g}^{2}}\right) +2k{_{n}}\left( \tau {_{g}^{2}-\tau {_{g}^{^{\prime }}}}%
\right) \right] + \\ 
2{\tau {_{g}^{^{\prime }}}}^{2}+2\tau {_{g}^{4}+}\left[ \left(
k_{n}^{^{\prime }}+k_{g}^{^{\prime }}\right) +\tau {_{g}}\left( k_{n}-k{_{g}}%
\right) \right] ^{2}%
\end{array}%
}+\frac{d\varphi ^{\ast }}{ds^{\ast }}
\end{equation*}%
\begin{corollary}
Consider that $\alpha \left( s\right) ~$is a geodesic curve, then the
following equation holds,
\end{corollary}~Consider that $\alpha \left( s\right) ~$is a principal line,
then the following equation holds,

$i)~\kappa ^{\ast }=\frac{\sqrt{2~\left( k_{n}^{2}+\tau {_{g}^{2}}\right) }}{%
\left( k_{n}+\tau {_{g}}\right) },$

$ii)~\ \tau ^{\ast }=\frac{-1}{\sqrt{2}}\frac{\left( k_{n}+\tau {_{g}}%
\right) ^{2}+\left( k_{g}k_{n}^{^{\prime }}-~k_{n}k_{g}^{^{\prime }}\right) 
}{\left( k_{n}^{2}+k{_{g}^{2}}\right) +\left( k_{n}^{^{\prime }}+{\tau {%
_{g}^{^{\prime }}}}\right) ^{2}\left( k_{n}+k{_{g}}\right) ^{-2}},$

$iii)~k_{g}^{\ast }=\frac{\sqrt{2~\left( k_{n}^{2}+\tau {_{g}^{2}}\right) }}{%
\left( k_{n}+\tau {_{g}}\right) }\cos \varphi ^{\ast },$

$iv)~\ k_{n}^{\ast }=\frac{\sqrt{2~\left( k_{n}^{2}+\tau {_{g}^{2}}\right) }%
}{\left( k_{n}+\tau {_{g}}\right) }\sin \varphi ^{\ast },$

$v)~\ \tau _{g}^{\ast }=\frac{-1}{\sqrt{2}}\frac{\left( k_{n}+\tau {_{g}}%
\right) ^{2}+\left( k_{g}k_{n}^{^{\prime }}-~k_{n}k_{g}^{^{\prime }}\right) 
}{\left( k_{n}^{2}+k{_{g}^{2}}\right) +\left( k_{n}^{^{\prime }}+{\tau {%
_{g}^{^{\prime }}}}\right) ^{2}\left( k_{n}+k{_{g}}\right) ^{-2}}+\frac{%
d\varphi ^{\ast }}{ds^{\ast }}.$

\subsection{$\mathbf{Tgn}{-~}$Smarandache Curves}

\begin{definition}
Let $M$ ~be an oriented surface in $E^{3}$ and let consider the arc - length
parameter curve $\alpha =\alpha \left( s\right) $ lying fully on $M$. Denote
the Darboux frame of $\alpha \left( s\right) $ $\left\{ \mathbf{T}{,\mathbf{g%
},}\mathbf{n}\right\} .$
\end{definition} Let $M$ ~be an oriented surface in $E^{3}$ and let consider
the arc - length parameter curve $\alpha =\alpha \left( s\right) $ lying
fully on $M$. Denote the Darboux frame of $\alpha \left( s\right) $ $\left\{ 
\mathbf{T},\mathbf{g},\mathbf{n}\right\} .$

$\mathbf{Tgn}$- Smarandache curve can be defined by 
\begin{equation}
\beta \left( s^{\ast }\right) =\frac{1}{\sqrt{3}}\left( \mathbf{T}+{\mathbf{g%
}+}\mathbf{n}\right) .  \label{3.16}
\end{equation}

\indent Now, we can investigate Darboux invariants of $\ \mathbf{Tgn}$ $~-~$%
Smarandache curve according to $\alpha =\alpha \left( s\right) .$
Differentiating (3.16) with respect to $s$, we get 
\begin{equation}
\beta ^{^{\prime }}=\frac{d\beta }{ds^{\ast }}\frac{ds^{\ast }}{ds}=\frac{-1%
}{\sqrt{3}}\left( \left( k_{n}+{k_{g}}\right) \mathbf{T}{+}\left( {\tau {%
_{g}-}k_{g}}\right) {\mathbf{g}-\left( {\tau {_{g}+}k_{n}}\right) }\mathbf{n}%
\right) ,  \label{3.17}
\end{equation}%
and 
\begin{equation*}
\mathbf{T}^{\ast }\frac{ds^{\ast }}{ds}=\frac{-1}{\sqrt{3}}\left( \left(
k_{n}+{k_{g}}\right) \mathbf{T}{+}\left( {\tau {_{g}-}k_{g}}\right) {\mathbf{%
g}-\left( {\tau {_{g}+}k_{n}}\right) }\mathbf{n}\right)
\end{equation*}%
where 
\begin{equation}
\frac{ds^{\ast }}{ds}=\sqrt{\frac{\left( k_{n}+{k_{g}}\right) ^{2}+\left( {%
\tau {_{g}-}k_{g}}\right) ^{2}+\left( {\tau {_{g}+}k_{n}}\right) ^{2}}{3}}.
\label{3.18}
\end{equation}%
The tangent vector of curve $\beta ~$can be written as follow, 
\begin{equation}
\mathbf{T}^{\ast }=\frac{-1}{\sqrt{\left( k_{n}+{k_{g}}\right) ^{2}+\left( {%
\tau {_{g}-}k_{g}}\right) ^{2}+\left( {\tau {_{g}+}k_{n}}\right) ^{2}}}%
\left( \left( k_{n}+{k_{g}}\right) \mathbf{T}{+}\left( {\tau {_{g}-}k_{g}}%
\right) {\mathbf{g}-\left( {\tau {_{g}+}k_{n}}\right) }\mathbf{n}\right) .
\label{3.19}
\end{equation}%
Differentiating (3.19) with respect to $s$, we obtain 
\begin{equation}
\frac{d\mathbf{T}^{\ast }}{ds^{\ast }}\frac{ds^{\ast }}{ds}=\frac{-1}{\left(
\left( k_{n}+{k_{g}}\right) ^{2}+\left( {\tau {_{g}-}k_{g}}\right)
^{2}+\left( {\tau {_{g}+}k_{n}}\right) ^{2}\right) ^{\frac{3}{2}}}\left(
\delta _{1}\mathbf{T}{+\delta _{2}\mathbf{g}+\delta _{3}}\mathbf{n}\right)
\label{3.20}
\end{equation}%
where 
\begin{equation}
\left\{ 
\begin{array}{l}
\delta _{1}={\left( k_{n}+{k_{g}}\right) }^{2}\left[ k_{g}\left( \tau {_{g}-}%
k_{g}\right) -k_{n}\left( \tau {_{g}+}k_{n}\right) \right] {+\left( k_{n}+{%
k_{g}}\right) }\left[ \left( \tau {_{g}-}k_{g}\right) \left( {\tau {%
_{g}^{^{\prime }}}-}k_{g}^{^{\prime }}\right) +\left( \tau {_{g}+}%
k_{g}\right) \left( {\tau {_{g}^{^{\prime }}}+}k_{n}^{^{\prime }}\right) %
\right] + \\ 
\left( \left( \tau {_{g}-}k_{g}\right) ^{2}+\left( \tau {_{g}+}k_{n}\right)
^{2}\right) \left[ k_{g}\left( \tau {_{g}-}k_{g}\right) -\left(
k_{g}^{^{\prime }}{+}k_{n}^{^{\prime }}\right) -k_{n}\left( \tau {_{g}+}%
k_{n}\right) \right] \\ 
\\ 
\delta _{2}=\left( \tau {_{g}-}k_{g}\right) ^{2}\left[ -k_{g}{\left( k_{n}+{%
k_{g}}\right) }-\tau {_{g}}\left( \tau {_{g}+}k_{n}\right) \right] {+}\left(
\tau {_{g}-}k_{g}\right) \left[ \left( k_{g}+k_{n}\right) \left(
k_{g}^{^{\prime }}{+}k_{n}^{^{\prime }}\right) +\left( \tau {_{g}+}%
k_{n}\right) \left( {\tau {_{g}^{^{\prime }}}+}k_{n}^{^{\prime }}\right) %
\right] + \\ 
\left( \left( k_{g}+k_{n}\right) ^{2}+\left( \tau {_{g}+}k_{n}\right)
^{2}\right) \left[ -k_{g}\left( k_{g}+k_{n}\right) +\left( k_{g}^{^{\prime }}%
{-\tau {_{g}^{^{\prime }}}}\right) -\tau {_{g}}\left( \tau {_{g}+}%
k_{n}\right) \right] \\ 
\\ 
\delta _{3}={\left( \tau {_{g}}+{k_{n}}\right) }^{2}\left[ \tau {_{g}}\left(
k_{g}{-}\tau {_{g}}\right) -k_{n}\left( k_{g}+k_{n}\right) \right] {+\left(
\tau {_{g}}+{k_{n}}\right) }\left[ -\left( k_{g}+k_{n}\right) \left(
k_{g}^{^{\prime }}+k_{n}^{^{\prime }}\right) -\left( \tau {_{g}-}%
k_{g}\right) \left( {\tau {_{g}^{^{\prime }}}-}k_{g}^{^{\prime }}\right) %
\right] + \\ 
\left( \left( k_{g}+k_{n}\right) ^{2}+\left( \tau {_{g}-}k_{g}\right)
^{2}\right) \left[ \tau {_{g}}\left( k_{g}-\tau {_{g}}\right) +\left( {\tau {%
_{g}^{^{\prime }}}+}k_{n}^{^{\prime }}\right) -k{_{n}}\left( k{_{g}+}%
k_{n}\right) \right]%
\end{array}%
\right.  \notag
\end{equation}%
Substituting (3.18) in (3.20), we get 
\begin{equation*}
\overset{\cdot }{\mathbf{T}^{\ast }}=\frac{\sqrt{3}}{\left( \left( k_{n}+{%
k_{g}}\right) ^{2}+\left( {\tau {_{g}-}k_{g}}\right) ^{2}+\left( {\tau {_{g}+%
}k_{n}}\right) ^{2}\right) ^{2}}\left( \delta _{1}\mathbf{T}{+\delta _{2}%
\mathbf{g}+\delta _{3}}\mathbf{n}\right) .
\end{equation*}%
Then, the curvature and principal normal vector field of curve $\beta ~$are
respectively, 
\begin{equation*}
\kappa ^{\ast }=\left\Vert \overset{\cdot }{\mathbf{T}^{\ast }}\right\Vert =%
\frac{\sqrt{2~\left( \delta _{1}^{2}{+\delta _{2}^{2}+\delta _{3}^{2}}%
\right) }}{\left( \left( k_{n}+{k_{g}}\right) ^{2}+\left( {\tau {_{g}-}k_{g}}%
\right) ^{2}+\left( {\tau {_{g}+}k_{n}}\right) ^{2}\right) ^{2}}
\end{equation*}%
\noindent and 
\begin{equation*}
\mathbf{N}^{\ast }=\frac{1}{\sqrt{\left( \delta _{1}^{2}{+\delta
_{2}^{2}+\delta _{3}^{2}}\right) }}\left( \delta _{1}\mathbf{T}{+\delta _{2}%
\mathbf{g}+\delta _{3}}\mathbf{n}\right) .
\end{equation*}%
\qquad \qquad \qquad On the other hand, we express 
\begin{equation*}
\mathbf{B}^{\ast }=\mathbf{T}^{\ast }{\times }\mathbf{N}^{\ast }=\frac{1}{%
\sqrt{\left( k_{n}+{k_{g}}\right) ^{2}+\left( {\tau {_{g}-}k_{g}}\right)
^{2}+\left( {\tau {_{g}+}k_{n}}\right) ^{2}}\sqrt{\left( \delta _{1}^{2}{%
+\delta _{2}^{2}+\delta _{3}^{2}}\right) }}\left\vert 
\begin{array}{ccc}
\mathbf{T} & {\mathbf{g}} & \mathbf{n} \\ 
-{\left( k_{n}+{k_{g}}\right) } & {k_{g}-}\tau {_{g}} & \tau {_{g}+}k_{n} \\ 
\delta _{1} & \delta {_{2}} & \delta {_{3}}%
\end{array}%
\right\vert .
\end{equation*}%
\qquad So, the binormal vector of curve $\beta ~$is

\begin{equation*}
\mathbf{B}^{\ast }=\frac{1}{\sqrt{\left( k_{n}+{k_{g}}\right) ^{2}+\left( {%
\tau {_{g}-}k_{g}}\right) ^{2}+\left( {\tau {_{g}+}k_{n}}\right) ^{2}}\sqrt{%
\lambda _{1}^{2}{+\lambda _{2}^{2}+\lambda _{3}^{2}}}}\left( \sigma _{1}%
\mathbf{T}{+\sigma _{2}\mathbf{g}+\sigma _{3}}\mathbf{n}\right)
\end{equation*}

where 
\begin{equation*}
\left\{ 
\begin{array}{l}
\sigma _{1}=\left( {k_{g}-}\tau {_{g}}\right) \delta {_{3}-}\left( \tau {%
_{g}+}k_{n}\right) \delta {_{2}} \\ 
\sigma {_{2}}=\left( \tau {_{g}+}k_{n}\right) {\delta _{1}+\left( k_{n}+{%
k_{g}}\right) }\delta {_{3}} \\ 
\sigma {_{3}}=-{\left( k_{n}+{k_{g}}\right) \delta {_{2}}-}\left( {k_{g}-}%
\tau {_{g}}\right) \delta _{1}.%
\end{array}%
\right.
\end{equation*}%
We differentiate (3.17) with respect to s in order to calculate the torsion

\begin{equation*}
\beta ^{^{\prime \prime }}=\frac{-1}{\sqrt{3}}\left\{ 
\begin{array}{l}
\left( k_{g}^{^{\prime }}+k_{n}^{^{\prime }}+k{_{g}}\left( {k_{g}-}\tau {_{g}%
}\right) +k_{n}\left( {\tau {_{g}+}k_{n}}\right) \right) \mathbf{T}{+}\left(
k_{g}{\left( k_{n}+{k_{g}}\right) }+\left( {\tau {_{g}^{^{\prime }}-}}%
k_{g}^{^{\prime }}\right) {+\tau {_{g}}}\left( {\tau {_{g}+}k_{n}}\right)
\right) {\mathbf{g}+} \\ 
\left( k_{n}{\left( k_{n}+{k_{g}}\right) }+{\tau {_{g}}}\left( {\tau {_{g}-}%
k_{g}}\right) -\left( {\tau {_{g}^{^{\prime }}}+}k_{n}^{^{\prime }}\right)
\right) \mathbf{n}%
\end{array}%
\right\} ,
\end{equation*}%
and similarly 
\begin{equation*}
\beta ^{^{\prime \prime \prime }}=\frac{-1}{\sqrt{3}}\left( \xi _{1}\mathbf{T%
}{+\xi _{2}\mathbf{g}+\xi _{3}}\mathbf{n}\right) ,
\end{equation*}%
where 
\begin{equation*}
\left\{ 
\begin{array}{l}
\xi _{1}=k_{n}^{^{\prime \prime }}+k_{g}^{^{\prime \prime
}}-2k_{g}^{^{\prime }}\left( k_{g}^{^{\prime }}-{\tau {_{g}^{^{\prime }}}}%
\right) -\left( k_{n}+k{_{g}}\right) \left( k_{n}^{2}+k_{g}^{2}\right)
+\left( {k_{g}-}\tau {_{g}}\right) \left( k_{g}^{^{\prime }}+k_{n}\tau {_{g}}%
\right) + \\ 
2k{_{n}}\left( k_{n}^{^{\prime }}+{\tau {_{g}^{^{\prime }}}}\right) +\left(
k_{n}+{\tau {_{g}}}\right) \left( k_{n}^{^{\prime }}-k_{g}\tau {_{g}}\right)
\\ 
\\ 
\xi _{2}={\tau {_{g}^{^{\prime \prime }}-}}k_{g}^{^{\prime \prime }}+2{\tau {%
_{g}}}\left( k_{n}^{^{\prime }}+{\tau {_{g}^{^{\prime }}}}\right) +2k{_{g}}%
\left( k_{n}^{^{\prime }}+k_{g}^{^{\prime }}\right) +\left( k_{n}+k{_{g}}%
\right) \left( k_{g}^{^{\prime }}-k_{n}\tau {_{g}}\right) + \\ 
\left( {\tau {_{g}+}k_{n}}\right) \left( k_{n}k{_{g}}+{\tau {_{g}^{^{\prime
}}}}\right) +\left( {k_{g}-}\tau {_{g}}\right) \left( {k_{g}^{2}+}\tau {%
_{g}^{2}}\right) \\ 
\\ 
\xi _{3}=-\left( {\tau {_{g}^{^{\prime \prime }}+}}k_{n}^{^{\prime \prime
}}\right) +2k{_{n}}\left( k_{n}^{^{\prime }}+k_{g}^{^{\prime }}\right)
+\left( {\tau {_{g}-}k_{g}}\right) \left( {\tau {_{g}^{^{\prime }}}}-k_{n}k{%
_{g}}\right) +\left( k_{n}+k{_{g}}\right) \left( k_{n}^{^{\prime
}}+k_{g}\tau {_{g}}\right) + \\ 
2\tau {_{g}}\left( {\tau {_{g}^{^{\prime }}}}-k_{g}^{^{\prime }}\right)
+\left( k_{n}+{\tau {_{g}}}\right) \left( k_{n}^{2}+\tau {_{g}^{2}}\right) .%
\end{array}%
\right.
\end{equation*}%
The torsion of curve $\beta ~$is

\begin{equation*}
\tau ^{\ast }=\frac{-1}{\sqrt{3}}\frac{%
\begin{array}{l}
\left( k_{n}+k{_{g}}\right) \left( k{_{g}}-\tau {_{g}}\right) \left( \tau {%
_{g}\xi }_{2}-k_{n}{\xi }_{1}\right) +\left( k{_{g}}-\tau {_{g}}\right)
^{2}\left( \tau {_{g}\xi }_{1}+k_{g}{\xi }_{3}\right) + \\ 
\left( k_{n}^{^{\prime }}+{\tau {_{g}^{^{\prime }}}}\right) \left[ \left( k{%
_{g}}-\tau {_{g}}\right) {\xi }_{1}+\left( k_{n}+k{_{g}}\right) {\xi }_{2}%
\right] - \\ 
\left( k_{n}^{^{\prime }}+k_{g}^{^{\prime }}\right) \left[ \left(
k_{n}^{^{\prime }}+\tau {_{g}}\right) {\xi }_{1}+\left( k{_{g}}-\tau {_{g}}%
\right) {\xi }_{3}-\left( k_{n}+\tau {_{g}}\right) {\xi }_{2}\right] + \\ 
\left( k{_{n}}+\tau {_{g}}\right) ^{2}\left( \tau {_{g}\xi }_{1}-k_{n}{\xi }%
_{2}\right) + \\ 
\left( k_{n}+\tau {_{g}}\right) \left( k_{n}+k{_{g}}\right) \left( \tau {%
_{g}\xi }_{3}+k_{g}{\xi }_{1}\right) +\left( k_{n}+\tau {_{g}}\right)
^{2}\left( k{_{g}\xi }_{3}-k_{n}{\xi }_{2}\right) + \\ 
\left( k_{g}-\tau {_{g}}\right) \left( k_{n}+\tau {_{g}}\right) \left( k{%
_{n}\xi }_{3}-k_{g}{\xi }_{2}\right) +\left( k_{n}+k{_{g}}\right) \left( {%
\tau {_{g}^{^{\prime }}}}-k_{g}^{^{\prime }}\right)%
\end{array}%
}{%
\begin{array}{l}
\left( k{_{g}}-\tau {_{g}}\right) ^{2}\left( k_{g}^{2}+\tau _{g}^{2}\right)
+2\left( k_{n}^{^{\prime }}+k_{g}^{^{\prime }}\right) \left[ k{_{g}}\left( k{%
_{g}}-\tau {_{g}}\right) +k{_{n}}\left( k_{n}+\tau {_{g}}\right) \right] +
\\ 
\left( k_{n}+\tau {_{g}}\right) ^{2}\left( k_{n}^{2}+\tau _{g}^{2}\right)
+2\left( k_{n}+\tau {_{g}}\right) \left[ \tau {_{g}}\left( {\tau {%
_{g}^{^{\prime }}}}-k_{g}^{^{\prime }}\right) +k{_{g}^{2}}\left( k_{n}+\tau {%
_{g}}\right) \right] + \\ 
\left( k_{n}+k{_{g}}\right) ^{2}\left( k_{n}^{2}+\tau _{g}^{2}\right)
+2\left( k_{n}+k{_{g}}\right) \left[ k{_{g}}\left( {\tau {_{g}^{^{\prime }}}}%
-k_{g}^{^{\prime }}\right) -k{_{n}}\left( {\tau {_{g}^{^{\prime }}}}%
+k_{g}^{^{\prime }}\right) \right] + \\ 
k_{n}\tau {_{g}}\left( \tau {_{g}}-k{_{g}}\right) +\left( k_{n}^{^{\prime
}}+k_{g}^{^{\prime }}\right) ^{2}+\left( k_{n}^{2}+\tau _{g}^{2}\right)
+\left( {\tau {_{g}^{^{\prime }}}}-k_{g}^{^{\prime }}\right) ^{2}+\left( {%
\tau {_{g}^{^{\prime }}}}-k_{n}^{^{\prime }}\right) ^{2}+ \\ 
-2\tau {_{g}}\left( \tau {_{g}}-k{_{g}}\right) \left( {\tau {_{g}^{^{\prime
}}}}+k_{n}^{^{\prime }}\right) .%
\end{array}%
}.
\end{equation*}%
The unit normal vector of surface $M$ and unit vector $\mathbf{g~}$of $\beta 
$ are as follow. Then, from (2.1) we obtain%
\begin{equation*}
\mathbf{g}^{\ast }=\frac{1}{\sqrt{\Phi }\Lambda }\left\{ \left( \sqrt{%
\Lambda }\cos \varphi ^{\ast }{\delta }_{1}+\sin \varphi ^{\ast }\sigma
_{1}\right) \mathbf{T}{+}\left( \sqrt{\Lambda }\cos \varphi ^{\ast }{\delta }%
_{2}+\sin \varphi ^{\ast }\sigma _{2}\right) \mathbf{g}{+}\left( \sqrt{%
\Lambda }\cos \varphi ^{\ast }{\delta }_{3}+\sin \varphi ^{\ast }\sigma
_{3}\right) \mathbf{n}\right\} ,
\end{equation*}%
and%
\begin{equation*}
\mathbf{n}^{\ast }=\frac{1}{\sqrt{\Phi }\Lambda }\left\{ \left( \sigma
_{1}\cos \varphi ^{\ast }-\sqrt{\Lambda }\sin \varphi ^{\ast }{\delta }%
_{1}\right) \mathbf{T}{+}\left( \sigma _{2}\cos \varphi ^{\ast }-\sqrt{%
\Lambda }\sin \varphi ^{\ast }{\delta }_{2}\right) \mathbf{g}{+}\left(
\sigma _{3}\cos \varphi ^{\ast }-\sqrt{\Lambda }\sin \varphi ^{\ast }{\delta 
}_{3}\right) \mathbf{n}\right\} .
\end{equation*}%
where $\Phi =\sqrt{\delta _{1}^{2}{+\delta _{2}^{2}+\delta _{3}^{2}}}%
,~\Lambda =\left( k_{n}+{k_{g}}\right) ^{2}+\left( {\tau {_{g}-}k_{g}}%
\right) ^{2}+\left( {\tau {_{g}+}k_{n}}\right) ^{2}$ and $\varphi ^{\ast }~$%
is the angle between the vectors ${\mathbf{g}}^{\ast }~$and $\mathbf{n}%
^{\ast }.~$Now, we can calculate geodesic curvature, normal curvature,
geodesic torsion of curve , so from (2.3) we get 
\begin{equation*}
k_{g}^{\ast }=\frac{\sqrt{2\left( \delta _{1}^{2}{+\delta _{2}^{2}+\delta
_{3}^{2}}\right) }}{\left( \left( k_{n}+{k_{g}}\right) ^{2}+\left( {\tau {%
_{g}-}k_{g}}\right) ^{2}+\left( {\tau {_{g}+}k_{n}}\right) ^{2}\right) ^{2}}%
\cos \varphi ^{\ast },
\end{equation*}%
and 
\begin{equation*}
k_{n}^{\ast }=\frac{\sqrt{2\left( \delta _{1}^{2}{+\delta _{2}^{2}+\delta
_{3}^{2}}\right) }}{\left( \left( k_{n}+{k_{g}}\right) ^{2}+\left( {\tau {%
_{g}-}k_{g}}\right) ^{2}+\left( {\tau {_{g}+}k_{n}}\right) ^{2}\right) ^{2}}%
\sin \varphi ^{\ast },
\end{equation*}%
and 
\begin{equation*}
\tau _{g}^{\ast }=\frac{-1}{\sqrt{3}}\frac{%
\begin{array}{l}
\left( k_{n}+k{_{g}}\right) \left( k{_{g}}-\tau {_{g}}\right) \left( \tau {%
_{g}\xi }_{2}-k_{n}{\xi }_{1}\right) +\left( k{_{g}}-\tau {_{g}}\right)
^{2}\left( \tau {_{g}\xi }_{1}+k_{g}{\xi }_{3}\right) + \\ 
\left( k_{n}^{^{\prime }}+{\tau {_{g}^{^{\prime }}}}\right) \left[ \left( k{%
_{g}}-\tau {_{g}}\right) {\xi }_{1}+\left( k_{n}+k{_{g}}\right) {\xi }_{2}%
\right] - \\ 
\left( k_{n}^{^{\prime }}+k_{g}^{^{\prime }}\right) \left[ \left(
k_{n}^{^{\prime }}+\tau {_{g}}\right) {\xi }_{1}+\left( k{_{g}}-\tau {_{g}}%
\right) {\xi }_{3}-\left( k_{n}+\tau {_{g}}\right) {\xi }_{2}\right] + \\ 
\left( k{_{n}}+\tau {_{g}}\right) ^{2}\left( \tau {_{g}\xi }_{1}-k_{n}{\xi }%
_{2}\right) + \\ 
\left( k_{n}+\tau {_{g}}\right) \left( k_{n}+k{_{g}}\right) \left( \tau {%
_{g}\xi }_{3}+k_{g}{\xi }_{1}\right) +\left( k_{n}+\tau {_{g}}\right)
^{2}\left( k{_{g}\xi }_{3}-k_{n}{\xi }_{2}\right) + \\ 
\left( k_{g}-\tau {_{g}}\right) \left( k_{n}+\tau {_{g}}\right) \left( k{%
_{n}\xi }_{3}-k_{g}{\xi }_{2}\right) +\left( k_{n}+k{_{g}}\right) \left( {%
\tau {_{g}^{^{\prime }}}}-k_{g}^{^{\prime }}\right)%
\end{array}%
}{%
\begin{array}{l}
\left( k{_{g}}-\tau {_{g}}\right) ^{2}\left( k_{g}^{2}+\tau _{g}^{2}\right)
+2\left( k_{n}^{^{\prime }}+k_{g}^{^{\prime }}\right) \left[ k{_{g}}\left( k{%
_{g}}-\tau {_{g}}\right) +k{_{n}}\left( k_{n}+\tau {_{g}}\right) \right] +
\\ 
\left( k_{n}+\tau {_{g}}\right) ^{2}\left( k_{n}^{2}+\tau _{g}^{2}\right)
+2\left( k_{n}+\tau {_{g}}\right) \left[ \tau {_{g}}\left( {\tau {%
_{g}^{^{\prime }}}}-k_{g}^{^{\prime }}\right) +k{_{g}^{2}}\left( k_{n}+\tau {%
_{g}}\right) \right] + \\ 
\left( k_{n}+k{_{g}}\right) ^{2}\left( k_{n}^{2}+\tau _{g}^{2}\right)
+2\left( k_{n}+k{_{g}}\right) \left[ k{_{g}}\left( {\tau {_{g}^{^{\prime }}}}%
-k_{g}^{^{\prime }}\right) -k{_{n}}\left( {\tau {_{g}^{^{\prime }}}}%
+k_{g}^{^{\prime }}\right) \right] + \\ 
k_{n}\tau {_{g}}\left( \tau {_{g}}-k{_{g}}\right) +\left( k_{n}^{^{\prime
}}+k_{g}^{^{\prime }}\right) ^{2}+\left( k_{n}^{2}+\tau _{g}^{2}\right)
+\left( {\tau {_{g}^{^{\prime }}}}-k_{g}^{^{\prime }}\right) ^{2}+\left( {%
\tau {_{g}^{^{\prime }}}}-k_{n}^{^{\prime }}\right) ^{2}+ \\ 
-2\tau {_{g}}\left( \tau {_{g}}-k{_{g}}\right) \left( {\tau {_{g}^{^{\prime
}}}}+k_{n}^{^{\prime }}\right) .%
\end{array}%
}+\frac{d\varphi ^{\ast }}{ds^{\ast }}.
\end{equation*}

\end{document}